\documentclass[letterpaper, 10 pt, conference]{ieeeconf}  
\IEEEoverridecommandlockouts                              

\usepackage{graphics}
\usepackage{graphicx}
\usepackage{epsfig}
\usepackage{amsmath,amsfonts,amssymb}
\usepackage{array}
\usepackage{multirow}
\usepackage{algorithm,float}
\usepackage{algorithmic}
\usepackage{mathrsfs}
\usepackage{epsfig,epstopdf}
\usepackage{textcomp}
\usepackage{booktabs}
\usepackage[colorlinks,linkcolor=blue,citecolor=blue]{hyperref}

\usepackage{hyperref}
\usepackage{authblk}
\usepackage{bm}
\usepackage{bbm,dsfont}
\usepackage{color}
\usepackage{pifont}
\usepackage{makecell}
\usepackage{url}
\usepackage{subfigure}
\usepackage{xr}
\usepackage{multirow}
\usepackage{indentfirst}
 \usepackage{cite}

\newtheorem{theorem}{Theorem}
\newtheorem{definition}{Definition}
\newtheorem{lemma}{Lemma}
\newtheorem{assumption}{Assumption}
\newtheorem{remark}{Remark}
\newtheorem{corollary}{Corollary}

\allowdisplaybreaks[4]
\newcommand{\bsx}{\boldsymbol{x}}
\newcommand{\bss}{\boldsymbol{s}}
\newcommand{\bse}{\boldsymbol{e}}

\newcommand{\bsX}{\boldsymbol{X}}
\newcommand{\bsy}{\boldsymbol{y}}
\newcommand{\sumn}{\sum\limits_{i=1}^n}

\newcommand{\sumT}{\sum\limits_{t=1}^T}

\title{\LARGE \bf
Quantized Distributed Online Projection-free Convex Optimization
}

\author{Wentao Zhang, Yang Shi, Baoyong Zhang, Kaihong Lu, Deming Yuan 
\thanks{
\emph{Corresponding author: Deming Yuan.}}
\thanks{W. Zhang, B. Zhang, D. Yuan are with School of Automation,  Nanjing University of Science and Technology,
        Nanjing 210094, Jiangsu, P. R. China (e-mail: iswt.zhang@gmail.com, baoyongzhang@njust.edu.cn, dmyuan1012@gmail.com).}%
\thanks{Y. Shi is with the Department of Mechanical Engineering, University of Victoria, Victoria, BC V8W 2Y2, Canada (e-mail: yshi@uvic.ca).}%
\thanks{K. Lu is with the College of Electrical Engineering and Automation,
Shandong University of Science and Technology, Qingdao 266590, China
(e-mail: khong\_lu@163.com)}}%

\begin{document}

\maketitle

\begin{abstract}
This paper considers online distributed convex constrained optimization over a time-varying multi-agent network. Agents in this network cooperate to minimize the global objective function through information exchange with their neighbors and local computation. Since the capacity or bandwidth of communication channels often is limited, a random quantizer is introduced to reduce the transmission bits. Through incorporating  this quantizer, we develop a quantized distributed online projection-free optimization algorithm, which can achieve the saving of  communication resources and computational costs. For different parameter settings of the quantizer, we establish the corresponding dynamic regret upper bounds of the proposed algorithm and reveal the trade-off between the convergence performance and the quantization effect. Finally, the theoretical results are illustrated by the simulation of distributed online linear regression problem.

\end{abstract}

\section{INTRODUCTION}
In recent years, online distributed convex optimization  has received ever-increasing attention from researchers because of its wide applications in many areas, such as machine learning, sensor networks,  smart grids,  etc.; see, e.g., \cite{yi2020distributed,shahrampour2017distributed, yuan2022distributedauto, yang2019survey, li2022survey, 8715380}. In such an online  optimization problem with constraint sets, various algorithms with projection operations have been developed, such as distributed online gradient descent \cite{sundhar2010distributed,6311406}. However, for some  high-dimensional and complex constrained optimization scenarios including multiclass classification \cite{pmlr-v70-zhang17g} and matrix completion \cite{hazan2012projection,7883821}, projection operations incur a heavy computational burden. On the contrary, projection-free algorithms have impressive advantages essentially due to the use of a linear oracle.

In \cite{pmlr-v70-zhang17g}, Zhang \emph{et al.}  earlier proposed an  online distributed projection-free algorithm and established the static regret upper bound  as $\mathcal{O} (T^{3/4})$. The works \cite{ pmlr-v119-wan20b, wan2021projection, thuang2022stochastic} further analyzed the static regret of some variants based on projection-free methods. In \cite{zhang2023dynamic} and \cite{10025380}, the  dynamic regret bounds were studied in distributed online projection-free algorithms under convex and nonconvex conditions, respectively.  Dynamic regret is a more stringent metric than static regret due to its dynamic reference sequence. However, the communication  channels between agents in \cite{zhang2023dynamic}, \cite{10025380}   are assumed to be perfect. In most applications,  the communication channels often have limited bandwidth or capacity, especially for the cases with scarce communication bandwidth or capacity \cite{cao2023communication}.

It is worth mentioning that quantized communication as a communication pattern
 can effectively reduce the number of communicated bits   to achieve the saving of communication resources \cite{cao2023communication}  through transmitting the information quantized by a dedicated quantizer. Currently, some distributed optimization algorithms with quantized communication have been developed \cite{xiong2022distributed, yi2014quantized, huang2016quantized, pu2016quantization, yuan2022distributedauto, 9224135,7891027, 9157925,9226092}, etc.  In \cite{xiong2022distributed}, Xiong \emph{et al.}  investigated the quantization effects on the convergence performance of the distributed quantized mirror descent algorithm. The works \cite{yi2014quantized, huang2016quantized, pu2016quantization} analyzed the quantized distributed off-line optimization algorithms based on subgradient and inexact proximal-gradient methods, respectively.
Doan \emph{et al.} \cite{9224135} considered a distributed off-line two-time-scale stochastic approximation algorithm under random quantization and established the almost sure convergence to the optimal solution for both convex and strongly convex loss function.
In \cite{7891027}, Li \emph{et al.} investigated the quantized distributed subgradient optimization algorithm with the dynamic encoding and decoding frameworks and proved that consensus optimization could be achieved under some mild conditions.
Further, for the online distributed optimization problem,  Yuan \emph{et al.} \cite{yuan2022distributedauto} proposed a distributed online bandit algorithm under quantized communication and established the static regret. Up to now, there are few research results considering distributed online optimization scenarios under quantized communication. The above analysis of the related literature and the state of the art motivates us to investigate the dynamic regret of distributed online projection-free algorithm under quantized communication, and the effect of quantizer parameters on regret bounds.

The main contributions of this work are two-fold. Firstly,  motivated by \cite{zhang2023dynamic} and \cite{yuan2022distributedauto},
 we develop a quantized distributed online projection-free optimization (Q-DOPFO) algorithm  for solving the distributed online constrained optimization problem over a  multi-agent network.  Meanwhile, the proposed algorithm  saves the communication resources and computational costs as compared to the algorithms with real-valued data and projection operations, respectively. Secondly, for different parameter settings of the quantizer, we establish the corresponding dynamic regret upper bounds of the proposed algorithm and reveal the trade-off between the convergence performance and the quantization effect. In particular, when the knowledge of $H_T$ is known, the optimal bound $\mathcal{O}( \sqrt{ T (1+H_T)} +D_T) $ can be achieved under proper parameter settings, where $T, H_T$ and $D_T$ represent the total time, function variation and gradient variation, respectively.

The remainder of the paper is organized as follows.  The problem statement, quantizer description, and necessary assumptions  are presented in Section II. Section III shows the algorithm design and convergence results. Sections IV and V give simulation examples and conclusion, respectively.


\textbf{Notation}:  ${\mathbb{R}}^n$  represents  the Euclidean space with $n$ dimensions. $[T]$ denotes $\{1,2,\ldots, T\}$.  $\|\boldsymbol{z}\|$ denotes the Euclidean  norm  of a vector $\boldsymbol{z}$. $ \lceil\cdot \rceil$ represents the round up function. The boundary of a set $\boldsymbol{X}$ is denoted as $\partial \boldsymbol{X}$. The element in the $i$-th row and $j$-th column of matrix $W$ is denoted as $[W]_{ij}$.  $[\boldsymbol{w}]_{i}$ denotes the $i$-th element of vector $\boldsymbol{w}$.  $\mathbb{B}_R^d:=\{\boldsymbol{z}\in\mathbb{R}^d |\  \| \boldsymbol{z}\| \leq R \}$ is the closed
Euclidean ball with a center point of origin and a radius of
$R$.
\section{Problem Formulation}
\subsection{The Optimization Problem}

  Consider a directed time-varying network $\mathcal{G}_t=\{\mathcal{V},\mathcal{E}_t,W_t \}$ that consists of $n$ agents, where $\mathcal{V} : = \{1,\ldots,n\} $, $ \mathcal{E}_t \subseteq \mathcal{V} \times \mathcal{V} $ denotes the edge set. In the network, agent $i$ can receive the information from the agents in its neighbor sets $\mathcal{N}_{i}^{in}(t)=\{j \mid(j, i) \in \mathcal{E}_t\} $.  $W_t \in \mathbb{R}^{n \times n}$ denotes the weighted  matrix and satisfies double stochasticity, i.e., $\sum_{j=1}^{n}  [W_t]_{ij}=\sum_{i=1}^{n} [W_t]_{ij}=1$,  $ \forall t\in [T], \forall i,j \in \mathcal{V }$, where $[W_t]_{ii}= 1-\sum_{j\in \mathcal{N}_{i}^{in}(t)} [W]_{ij}$. There exists a constant $ \zeta >0 $ such that $ [W_t]_{ij}>\zeta, t \in [T]$ holds when $j \in \mathcal{N}_{i}^{in}(t)\cup \{ i\}$, and $[W_t]_{ij}=0$ otherwise.
   The distributed online optimization problem is described as follows:
\begin{flalign}\label{problem definition}
\min\limits_{\boldsymbol{x}_t \in \boldsymbol{X}}\,\,\sum\limits_{t=1}^T F_t (\boldsymbol{x}_t)
\end{flalign}
where $F_t (\boldsymbol{x})= \sum_{i=1}^n{f_{i,t}}( \boldsymbol{x} )$, the function $f_{i,t}$ is convex over   the convex and compact set $\bsX \in$ ${\mathbb{R}}^d$. Agents in the network cooperate to search for the global optima of Problem (\ref{problem definition}) through local computation and information exchange with neighbor agents.
Generally, the metric  $\textbf{Regret}_d^j (T)$ defined in (\ref{Regret-j})  is used to measure the algorithm performance, which represents  the difference between the cumulative cost $F_t { (\boldsymbol{x}_{j,t})}$  of the agent $j$ over time $T$ and the cumulative cost at benchmark sequence $\bsx_{t}^* \in \bsX$.
\begin{flalign}\label{Regret-j}
\textbf{Regret}_d^j (T) =   \sum_{t=1}^T  F_t { (\boldsymbol{x}_{j,t})} -  \sum\limits_{t=1}^T F_t  (\bsx_t^*)
\end{flalign}
where  $\bsx_{t}^* \in {\arg\min}_ {\boldsymbol{x} \in \boldsymbol{X}}  F_t(\boldsymbol{x})$. Due to this varying benchmark $\bsx_t^*$, dynamic regret is more stringent than static regret and has wider application scenarios, such as target tracking. It is well known that  the upper bound of (\ref{Regret-j})  generally depends on the regularity of the optimization problem. Considering this fact, we define the following \emph{function variation} $H_T$ and \emph{gradient variation} $D_T$.
\begin{flalign} \label {H_TD_T}
H_T&:= \sum_{t=1}^{T-1} f_t^{sup},\
D_T:= \sum_{t=1}^{T-1} g_t^{sup}
\end{flalign}
where $f_t^{sup} =\max_{i \in \mathcal{V}} \max_{\bsx \in \boldsymbol{X}} |  f_{i,t+1}(\bsx)- f_{i,t}(\bsx)|$, $g_t^{sup}=\max_{i \in \mathcal{V}} \max_{\bsx \in \boldsymbol{X}}\| \nabla f_{i,t+1}(\bsx)-  \nabla f_{i,t}(\bsx)\| $.

 Our objective is to design a distributed online algorithm with quantized communication for Problem  (\ref{problem definition}) that achieves sublinear dynamic regret of every agent $j \in \mathcal{V} $.
\subsection{Random Quantizer}
In this section, the following random quantizer is introduced to  ensure that each agent in the network uses its quantized information to communicate with its neighbors.

\begin{definition} [\cite{yuan2022distributedauto}] \label{quantizer definition}
$\mathds{Q}_t (\bsy)\in \mathbb{R}^d$ is the time-varying random quantizer of a vector $\bsy \in \mathbb{R}^d$ if it satisfies that
\begin{align}
&\mathbb{E}[\mathds{Q}_t (\bsy)] = \bsy,\ \mathbb{E}[\|\mathds{Q}_t (\bsy)- \bsy \|^2] \leq \epsilon_{d,k_t} \|\bsy \|^2,  t \in [T]
\end{align}
where $\epsilon_{d,k_t}$ denotes a quantization resolution that  is dependent on the qunantization levels $k_t$ and the dimension $d$.
\end{definition}
\begin{remark} \label{remark 1}
Several common quantizers are naturally special cases of this random quantizer, such as randomized gossip \cite{koloskova2019decentralized}, rescaled unbiased estimators \cite{koloskova2019decentralized}, stochastic $k$-level quantization \cite{yuan2022distributedauto}, probabilistic quantizer \cite{xiong2022distributed}. We show the probabilistic quantizer in \cite{xiong2022distributed} as an example. Denote  $\mathds{Q}_t (\bsy) = [\mathds{Q}_t (a_1), \mathds{Q}_t (a_2), \ldots, \mathds{Q}_t (a_d)]^T $, where  $a_i=[\bsy]_i, i\in [d]$. Then, for $[\bsy]_i, i\in [d], t \in [T] $ , we have
\begin{align}
\mathds{Q}_t (a_i)=
\left\{\begin{array}{rcl}
{\overline{ a_i}}^t, \quad  w.p. \ (a_i - {\underline{ a_i}}^t ) k_t, \\
{ \underline{a_i}}^t, \quad w.p. \ (\overline{a_i}^t -a_i) k_t.
\end{array}\right.
\end{align}
where ${\overline{ a_i}}^t$ and ${\underline{a_i}}^t$ are the  round up and down $a_i$ to the nearest integer multiple of $1/ k_t$, respectively. It is not hard to note that the probabilistic quantizer satisfies Definition \ref{quantizer definition} with $\epsilon_{d,k_t}= d / (4 {k_t}^2)$.
\end{remark}
\begin{remark}
According to  Definition \ref{quantizer definition}, $\epsilon_{d,k_t}$ has a wide range of values and when its value is smaller, the quantized data is closer to the real-value data. Note that large values of $\epsilon_{d,k_t}$ are allowed at the early stages of the running algorithm, which means that the quantized data at this stage is coarser and less precise than the real-value data.
In order to  achieve the sublinear dynamic regret, a sublinearly convergent sequence $\{ \epsilon_{d,k_t} \}$ over time $t$ is desired and necessary, which
can be verified  in the following sections.
\end{remark}
\subsection{Some Assumptions}
Some necessary assumptions are needed to facilitate the following algorithm development.
\begin{assumption}\label{assump: network}
 The union  $\bigcup_{i=k Q+1}^{(k+1)Q} \mathcal{G}_{i}$  is strongly connected for some positive integer $Q$ and every integer $k \geq 0$.
\end{assumption}
\begin{assumption}\label{assump: ball}
The constraint set $\bsX \subset {\mathbb{R}}^d$ is convex and compact and satisfies that $\bsX \subseteq \mathbb{B}_R^d, R>0$.
\end{assumption}
\begin{assumption}\label{assump: lips funon}
 The function $f_{i,t}$ is $L_{X}$-Lipschitz,  i.e.,  $
 |f_{i,t}(\boldsymbol{x}_{1})-f_{i,t}(\boldsymbol{x}_{2})| \leq L_{X} \|\boldsymbol{x}_{1}-\boldsymbol{x}_{2}\|$, $\forall \boldsymbol{x}_1,\boldsymbol{x}_2 \in \boldsymbol{X}$,
where $L_{X}$  is a known positive constant.
\end{assumption}
\begin{assumption}\label{assump: lips grad}
The gradient $\nabla f_{i,t} (\boldsymbol{x})$ is $G_{X}$-Lipschitz, i.e., $
\| \nabla f_{i,t} (\boldsymbol{x}_1) - \nabla f_{i,t} (\boldsymbol{x}_2) \|\leq {G_X} \| \boldsymbol{x}_1-\boldsymbol{x}_2\|, \forall \boldsymbol{x}_1,\boldsymbol{x}_2 \in \boldsymbol{X}$, which is equivalent to $ f_{i,t}(\boldsymbol{x}_{1})-f_{i,t}(\boldsymbol{x}_{2})  \leq \langle \nabla f_{i,t} (\bsx_2), \bsx_1-\bsx_2\rangle
  \quad + \frac{G_X}{2} \|\boldsymbol{x}_{1}-\boldsymbol{x}_{2}\|^2$.
\end{assumption}
\begin{remark}
 Assumptions \ref{assump: network}-\ref{assump: lips funon} are common in the literature (see \cite{nedic2008distributed, yi2020distributed}, \cite{7883821,besbes2015non}, etc.) on centralized and distributed optimization. The purpose of assuming $\bsX \subseteq \mathbb{B}_R^d$ is to ensure that the variance of the random quantizer is bounded, i.e., $\mathbb{E}[\|\mathds{Q}_t (\bsy)- \bsy \|^2] \leq \epsilon_{d,k_t} R^2$, which is a necessary precondition.
It is worth noting that
 Assumption \ref{assump: lips funon} implies  $\|\nabla f_{i,t}(\bsx) \|\leq L_X$ according to Lemma 2.6 in \cite{shalev2011online}.
\end{remark}
\section{Algorithm Design and Convergence Analysis}
\subsection{Algorithm Q-DOPFO}
In this section, we  develop Algorithm Q-DOPFO, which is illustrated in Algorithm \ref{algorithm 1}. The key ingredients of the proposed algorithm include: 1) the quantized data $\mathds{Q}_t (\bsx_{j,t})$ and $\mathds{Q}_t[\nabla f_{i,t}  (\hat{\boldsymbol{x}}_{i,t})]$, instead of real-valued data, are utilized to perform consensus steps;
2) gradient tracking technique is introduced to correct the gradient change of loss function by using the global gradient estimation $\widehat{\bss} _{i,t}$ instead of individual agent gradients;
3) the  decision variable $\bsx_{i, t+1}$ is updated through a linear step.
It is worth noting that the use of the random quantizer and projection-free oracle in the proposed algorithm can effectively save  communication and computing resources of multi-agent systems.
\begin{algorithm}[H]
	\renewcommand{\algorithmicrequire}{\textbf{Initialize:} }
	\caption{ (Q-DOPFO) Quantized Distributed Online Projection-free Optimization }
	\label{algorithm 1}
	\begin{algorithmic}[1]
		\REQUIRE Initial  variables $\boldsymbol{x}_{i,1} \in \boldsymbol{X}$ and parameter $0 < \alpha \leq 1.$
		
		\FOR {$t=1,2,\cdots,T$}
             \FOR {Each agent $i \in \mathcal{V}$}
             \STATE Agent $i$ quantizes its state $\bsx_{i,t}$ and executes
             \IF{$\mathds{Q}_t(\bsx_{i,t}) \notin \boldsymbol{X}$}
            \STATE $\mathds{Q}_t(\bsx_{i,t}) = \bsx_{i,t}$.
            \ENDIF

             \STATE Agent $i$  receives the quantized data $\mathds{Q}_t (\bsx_{j,t})$ from its neighbors  $j \in \mathcal{N}_{i}^{in}(t)$, and  updates

                \begin{center}
                \quad $\boldsymbol{\hat{x}}_{i,t}= [W_t]_{ii} \mathds{Q}_t(\boldsymbol{x}_{i,t})+ \underset{j \in \mathcal{N}_{i}^{in}}{\sum}{[W_t]_{ij} \mathds{Q}_t(\boldsymbol{x}_{j,t})}.$
                           \setlength{\parskip}{0.4em}
                \end{center}

\STATE  After the gradient value $\nabla f_{i,t}  (\boldsymbol{\hat{x}}_{i,t})$  is revealed, agent $i$ obtains $\mathds{Q}_t[\nabla f_{i,t}  (\boldsymbol{\hat{x}}_{i,t})]$ and executes gradient tracking steps:
            \IF{$t=1 $}
            \STATE $\overline{\nabla} f_{i,1} =  \mathds{Q}_1[\nabla f_{i,1}  (\hat{\boldsymbol{x}}_{i,1})]$,
            \ELSE
            \STATE $\overline{\nabla} f_{i,t} = \widehat{\bss} _{i,t-1} + \mathds{Q}_t[\nabla f_{i,t}  (\hat{\boldsymbol{x}}_{i,t})]$\\
            $\quad \quad \quad \quad  \quad \quad \quad \quad \quad - \mathds{Q}_{t-1}[\nabla f_{i,t-1}  (\hat{\boldsymbol{x}}_{i,t-1})].$
            \ENDIF

               \begin{center}
               \quad $\widehat{\bss} _{i,t}= {[W_t]_{ii}} \overline{\nabla} f_{i,t}+ \sum_{j \in \mathcal{N}_{i}^{in}(t)}{[W_t]_{ij}} \overline{\nabla} f_{j,t},$
               \end{center}

		\STATE Frank-Wolfe step: update

             \begin{center}
               \quad $\boldsymbol{v}_{i,t} = \underset{\boldsymbol{x} \in \boldsymbol{X}}{\arg\min} \left<\boldsymbol{x}, \widehat{\bss} _{i,t} \right>, $

               \quad $\boldsymbol{x}_{i,t+1}= \hat{\boldsymbol{x}}_{i,t} +\alpha (\boldsymbol{v}_{i,t}- \hat{\boldsymbol{x}}_{i,t}).$
              \end{center}

		\ENDFOR
		\ENDFOR
	\end{algorithmic}
\end{algorithm}

 In some extreme situations, such as $\boldsymbol{x}_{i,t} \in \partial \boldsymbol{X}$ at time $t$, $\mathds{Q}_t(\bsx_{i,t})$ may occasionally violate the constraint set due to the quantizer. However, because of the variability of $x_t^*$ over the time and the randomness of the quantizer, $\mathds{Q}_t(\bsx_{i,t})$ usually does not always violate set $\bsX$. To ensure that the updated decipsion $\bsx_{i,t+1}$ is always feasible, we require that the quantized state $\mathds{Q}_t(\bsx_{i,t})$ is in set $\bsX$ for all $t$, i.e. the step 4 of Algorithm \ref{algorithm 1}.

\subsection{Main Convergence Results}
In this section, some lemmas and the bound of dynamic regret defined in (\ref{Regret-j}) for Algorithm \ref{algorithm 1} are  established. To facilitate the analysis, we define as follows the transition matrix $\Phi(t, s)=W_t W_{t-1} \ldots W_{s}$, for all $t, s \ \text{with} \  t \geq s \geq 1$, the running average vectors $\bsx_{a,t}$, $\boldsymbol{v}_{a,t}$, the quantization errors $\boldsymbol{e}_{i,t},\boldsymbol{\theta}_{i,t}$   and  the difference of quantized gradient   $\boldsymbol{\nabla}_{i,t}^Q$.
\begin{gather} \label{average-delta}
\left\{\begin{array}{rcl}
\boldsymbol{x}_{a,t}&=&\frac{1}{n} \sum_{i=1}^n \boldsymbol{x}_{i,t},
\boldsymbol{v}_{a,t}=\frac{1}{n} \sum_{i=1}^n \boldsymbol{v}_{i,t} \\
\boldsymbol{e}_{i,t}&= &\mathds{Q}_t(\bsx_{i,t})-\bsx_{i,t} \\ \boldsymbol{\theta}_{i,t}&=&\mathds{Q}_t[\nabla f_{i,t}  (\hat{\boldsymbol{x}}_{i,t})]-\nabla f_{i,t}  (\hat{\boldsymbol{x}}_{i,t})\\
\boldsymbol{\nabla}_{i,t}^Q &=& \mathds{Q}_t[\nabla f_{i,t}  (\hat{\boldsymbol{x}}_{i,t})] - \mathds{Q}_{t-1} [\nabla f_{i,t-1}  (\hat{\boldsymbol{x}}_{i,t-1})]
\end{array}\right.
\end{gather}
\begin{lemma} \label{consistency}
Let the decision sequence $\{ \bsx_{i,t}\}$ be generated by Algorithm \ref{algorithm 1}. Then, under Assumptions \ref{assump: network} and \ref{assump: ball},  we have for $T\geq 2$ that
\begin{flalign}
& \sumT \sumn \mathbb{E} [\| \hat{\bsx}_{i,t}- \bsx_{a,t}\|]  \leq    \frac{n\Gamma}{1-\sigma} \sum\limits_{j=1}^n  \| {\boldsymbol{x}}_{j,1} \| +\alpha T\frac{2n^2R\Gamma}{1-\sigma} \nonumber \\
  &+\left(1+\frac{n \Gamma \sigma}{1-\sigma}\right)\sum_{t=1}^T \sum_{i=1}^n  \mathbb{E} [\|\bse_{i,t}\|]
\end{flalign}
where  $\sigma=(1-{\zeta}/ { 4 n^{2}})^{1 /  Q}, \Gamma=(1-{\zeta} /{ 4 n^{2}})^{(1-2Q) / Q}$.
\end{lemma}
\begin{lemma} \label{grad diffience}
Let the sequence $\{ \widehat{\bss} _{i,t}, {\nabla} f_{i,t} (\hat{\boldsymbol{x}}_{i,t}) \}$ be generated by Algorithm \ref{algorithm 1}. Then, under Assumptions \ref{assump: network} and \ref{assump: lips grad}, we have for any $T\geq 2$ that
\begin{flalign} \label{condition-lem4}
& \sumT \sumn \mathbb{E}\left[\left\|  \widehat{\bss} _{i,t} - \frac{1}{n}  {\nabla} F_{t} ({\boldsymbol{x}}_{a,t})\right\|\right]  \\
& \leq  C_1   + G_X C_2 \sum\limits_{t=1}^{T} \sumn \mathbb{E} \left[\left\|  \hat{\bsx}_{i,t}-\bsx_{a,t}  \right\| \right] + nL_X C_2 \sumT  \sqrt{\epsilon_{d,k_t}} \nonumber \\
&+ \frac{  n \Gamma G_X}{1-\sigma}  \sumT \sumn \mathbb{E}[\|\boldsymbol{e}_{i,t} \|] + \frac{  n^2 \Gamma }{1-\sigma} D_T+ \frac{ 2 n^2 \Gamma R G_X}{1-\sigma} \alpha T \nonumber
\end{flalign}
where $C_1=\frac{ \sigma n \Gamma \sqrt{\epsilon_{d,k_1}} + n \Gamma }{1-\sigma}  \sum_{i=1}^n \| {\nabla} f_{i,1} (\hat{\boldsymbol{x}}_{i,1})\|,  C_2=\frac{ 2n \Gamma}{1-\sigma}+1$.
\end{lemma}

\begin{theorem}\label{theorem 1}
Let the decision sequence $\{ \bsx_{i,t}\}$ be generated by Algorithm \ref{algorithm 1} and suppose Assumptions \ref{assump: network}-\ref{assump: lips grad} hold. Then,     for $T\geq 2$ and $j \in \mathcal{V}$,  the regret is bounded as follows:
\begin{flalign}\label{Eq Theorem1}
\mathbb{E}[\textbf{Regret}_d^j (T)]  &\leq D_1  +D_2 \alpha T +D_3 \sumT \sqrt{\epsilon_{d,k_t}}+ \frac{D_4}{\alpha} \nonumber \\
& \quad + \frac{D_5}{\alpha}\sumT \epsilon_{d,k_t} +\frac{2n}{\alpha} H_T  +D_6 D_T
\end{flalign}
where
\begin{align*}
& D_1=n L_{X} \sum_{i=1}^n \|\boldsymbol{x}_{i,1}-\boldsymbol{x}_{a,1}\|
 +\frac{n\Gamma E_0}{1-\sigma} \sum_{i=1}^n  \| {\boldsymbol{x}}_{i,1} \|+ 4RC_1,  \\
 &D_2=4n  R(nL_X+G_X R) + \frac{ 2n^2 R  \Gamma E_0 }{1-\sigma} +\frac{8 n^2 \Gamma G_X R^2}{2},\\
 &D_3=nR E_0 \left( 1+\frac{n\Gamma \sigma}{1-\sigma} \right) +n^2L_X R +4nR L_X C_2 \\
 &+ \frac{4n^2 \Gamma G_X R^2}{1-\sigma}, D_4=2n L_X R,  D_5=nG_X R^2, \\ & D_6=\frac{4 n^2 R \Gamma}{1-\sigma}+ 2n R, E_0=4R C_2 G_X+n L_X. \\
\end{align*}

\end{theorem}
\noindent{\em Proof.}
Based on  Algorithm \ref{algorithm 1} and  double stochasticity of $W_{t-1}$, we obtain that
\begin{flalign}
{\boldsymbol{x}}_{a,t}&= \frac{1}{n} \sum_{i=1}^n [\hat{\boldsymbol{x}}_{i,t-1}+ \alpha ( {\boldsymbol{v}}_{i,t-1}-\hat{\boldsymbol{x}}_{i,t-1})] \nonumber \\
&=\frac{1-\alpha}{n} \sum_{i=1}^n \mathds{Q}_{t-1} (\boldsymbol{x}_{i,t-1}) +\alpha {\boldsymbol{v}}_{a,t-1} \nonumber \\
&= \frac{1-\alpha}{n} \sum_{i=1}^n \boldsymbol{e}_{i,t-1} +(1-\alpha) \bsx_{a,t-1}+\alpha {\boldsymbol{v}}_{a,t-1}.
\end{flalign}

Thus, according to Assumptions \ref{assump: ball} and \ref{assump: lips funon}, for any $t\geq2$, we have that
\begin{flalign}\label{Proof Theorem1- 2-a1}
&F_t {  (\boldsymbol{x}_{j,t})} - F_t { (\boldsymbol{x}_{a,t})} \nonumber \\
&\leq nL_{X}      \|\boldsymbol{x}_{j,t}-\boldsymbol{x}_{a,t}\| \nonumber \\
&\leq   n L_{X} \sum\limits_{i=1}^n   \|\boldsymbol{x}_{i,t}-\boldsymbol{x}_{a,t}\| \nonumber \\
&= n L_{X} \sum\limits_{i=1}^n  \|\hat{\boldsymbol{x}}_{i,t-1}-\boldsymbol{x}_{a,t-1}+ \alpha (\boldsymbol{v}_{i,t-1}-\hat{\boldsymbol{x}}_{i,t-1})\nonumber \\
&\quad- \alpha (\boldsymbol{v}_{a,t-1}-{\boldsymbol{x}}_{a,t-1})- \frac{1-\alpha}{n} \sum_{i=1}^n \boldsymbol{e}_{i,t-1}\|\nonumber \\
&\leq  n L_{X} \sum\limits_{i=1}^n  \|\hat{\boldsymbol{x}}_{i,t-1}-\boldsymbol{x}_{a,t-1}\|+ n L_{X} \sum_{i=1}^n \|\boldsymbol{e}_{i,t-1}\|\nonumber \\
&\quad+ 4n^2 L_{X}\alpha  R.
\end{flalign}

Recalling the regret notion defined in  (\ref{Regret-j}) and combining (\ref{Proof Theorem1- 2-a1}), we obtain that
\begin{flalign}\label{Proof Theorem1- 1-a2}
&\mathbb{E}[\textbf{Regret}_d^j (T)]\nonumber \\
&\leq \sum_{t=1}^T \mathbb{E}[ F_t { (\boldsymbol{x}_{j,t})} -F_t { (\boldsymbol{x}_{a,t})} ] +  \sum_{t=1}^T \mathbb{E}[F_t { (\boldsymbol{x}_{a,t})}- F_t  (\bsx_t^*)]\nonumber \\
&\leq   n L_{X}
 \sum\limits_{i=1}^n \|\boldsymbol{x}_{i,1}-\boldsymbol{x}_{a,1}\| +  n L_{X}\sum\limits_{t=1}^{T-1} \sum\limits_{i=1}^n \mathbb{E}[\|\hat{\boldsymbol{x}}_{i,t}-\boldsymbol{x}_{a,t}\|] \nonumber \\
  &\quad + n^2 L_{X} R \sumT \sqrt{\epsilon_{d, k_t}}
 +  \sum_{t=1}^T \mathbb{E}[F_t { (\boldsymbol{x}_{a,t})}- F_t  (\bsx_t^*)] \nonumber \\
 &\quad +4\alpha T n^2 L_{X} R
\end{flalign}
where the last inequality is obtained by using the fact $\mathbb{E} [\| \bse_{i,t}\|] \leq \sqrt{\mathbb{E} [\| \bse_{i,t}\|^2]}  \leq \sqrt{\epsilon_{d,k_t}  \| \bsx_{i,t}\|^2  }\leq R \sqrt{\epsilon_{d,k_t}}$.

Next, we aim to bound the term $\sum_{t=1}^T \mathbb{E}[F_t { (\boldsymbol{x}_{a,t})}- F_t  (\bsx_t^*)]$ in (\ref{Proof Theorem1- 1-a2}). By using Assumption \ref{assump: lips grad}, we have
\begin{flalign} \label{proof lem key 1}
&F_{t+1} (\boldsymbol{x}_{a,t+1}) -F_{t+1} (\boldsymbol{x}_{a,t})  \nonumber\\
&\leq  \left< \nabla F_{t+1} (\boldsymbol{x}_{a,t}), \boldsymbol{x}_{a,t+1}- \boldsymbol{x}_{a,t} \right> + \frac{nG_X}{2} \| \boldsymbol{x}_{a,t+1}- \boldsymbol{x}_{a,t} \|^2 \nonumber \\
&\leq  \left\langle \nabla F_{t+1} (\boldsymbol{x}_{a,t}), \frac{1-\alpha}{n} \sum_{i=1}^n \boldsymbol{e}_{i,t} +\alpha (\boldsymbol{v}_{a,t}- \boldsymbol{x}_{a,t}) \right\rangle \nonumber \\
&\quad + \frac{nG_X}{2} \left\| \frac{1-\alpha}{n} \sum_{i=1}^n \boldsymbol{e}_{i,t} +\alpha (\boldsymbol{v}_{a,t}- \boldsymbol{x}_{a,t}) \right\|^2 \nonumber \\
&\leq  \alpha \sum_{i=1}^n \left\langle \frac{1}{n}\nabla F_{t+1} (\boldsymbol{x}_{a,t}),  \boldsymbol{v}_{i,t}- \boldsymbol{x}_{a,t} \right\rangle + \frac{1-\alpha}{n} \sum_{i=1}^n \Upsilon_{i,t} \nonumber \\
&\quad + {nG_X} \left[ \frac{(1-\alpha)^2}{n^2} \left\| \sum_{i=1}^n \boldsymbol{e}_{i,t}\right\|^2 +\alpha^2 \left\|\boldsymbol{v}_{a,t}- \boldsymbol{x}_{a,t} \right\|^2 \right]  \nonumber \\
&\leq  \alpha \sum_{i=1}^n \left\langle \frac{1}{n}\nabla F_{t+1} (\boldsymbol{x}_{a,t}),  \boldsymbol{v}_{i,t}- \boldsymbol{x}_{a,t} \right\rangle + \frac{1-\alpha}{n} \sum_{i=1}^n \Upsilon_{i,t} \nonumber \\
&\quad + {G_X}  \sum_{i=1}^n \left\|\boldsymbol{e}_{i,t}\right\|^2 +4nG_X R^2 \alpha^2
\end{flalign}
where $ \Upsilon_{i,t} := \langle \nabla F_{t+1} (\boldsymbol{x}_{a,t}), \boldsymbol{e}_{i,t}\rangle $ and the last inequality is obtained by using the fact $\| \sum_{i=1}^n \boldsymbol{e}_{i,t}\|^2 \leq n\sum_{i=1}^n \|\boldsymbol{e}_{i,t}\|^2  $ \cite{9184135}.
It can be further verified that
\begin{flalign} \label{proof lem key 2}
& \left< \frac{1}{n}\nabla F_{t+1} (\boldsymbol{x}_{a,t}), \boldsymbol{v}_{i,t}- \boldsymbol{x}_{a,t} \right> \nonumber \\
&\leq  \left< \frac{1}{n}\nabla F_{t+1} (\boldsymbol{x}_{a,t})- \widehat{\bss} _{i,t}, \boldsymbol{v}_{i,t}- \boldsymbol{x}_{a,t} \right> \nonumber \\
&\quad + \left< \widehat{\bss} _{i,t}, \bsx_t^*- \boldsymbol{x}_{a,t} \right> \nonumber \\
&\leq \frac{ 2R}{n} \left\|  \nabla F_{t+1} (\bsx_{a,t}) - \nabla F_{t} (\bsx_{a,t}) \right\|+{ 4R} \left\| \frac{1}{n}\nabla F_{t} (\bsx_{a,t}) \right.\nonumber \\
&\quad \left. - \widehat{\bss} _{i,t} \right\| + \frac{1}{n} \left< \nabla F_{t} (\bsx_{a,t}), \bsx_t^*- \bsx_{a,t} \right>
   \nonumber \\
&\leq  { 2R} g_{t,sup} +{ 4R} \left\| \frac{1}{n}\nabla F_{t} (\bsx_{a,t}) - \widehat{\bss} _{i,t} \right\| \nonumber \\
&\quad   +  \frac{1}{n}\left[ F_{t}(\bsx_t^*) - F_{t} (\bsx_{a,t}) \right]
\end{flalign}
where the first inequality is obtained by utilizing the following  optimality condition:
$
\langle \widehat{\bss} _{i,t}, \bsx_t^* \rangle  \geq  \langle \widehat{\bss} _{i,t}, \boldsymbol{v}_{i,t} \rangle
$
and the last inequality is derived  based on  the convexity condition of $F_t( \bsx)$ together with   Assumption \ref{assump: lips grad}. Then, it follows from (\ref{proof lem key 1}) and  (\ref{proof lem key 2}) that
\begin{flalign} \label{proof lem key 3}
&F_{t+1} (\boldsymbol{x}_{a,t+1})-F_{t+1} (\boldsymbol{x}_{a,t}) \nonumber \\
&\leq  { 2 \alpha n R} g_{t,sup} + \alpha \left[ F_{t}(\bsx_t^*) - F_{t} (\bsx_{a,t}) \right] + \Omega_t
\end{flalign}
where $\Omega_t={ 4\alpha R} \sum_{i=1}^n \| \frac{1}{n}\nabla F_{t} (\bsx_{a,t}) - \widehat{\bss} _{i,t} \| + \frac{1-\alpha}{n} \sum_{i=1}^n \Upsilon_{i,t} + {G_X}  \sum_{i=1}^n \left\|\boldsymbol{e}_{i,t}\right\|^2 +4nG_X R^2 \alpha^2$.
Through simplifying (\ref{proof lem key 3}) by using the method similar to Lemma 4 in \cite{zhang2023dynamic}, we obtain that
\begin{flalign} \label{proof lem key 6-a6}
\alpha \sum\limits_{t=1}^{T} \mathbb{E} \left[ F_{t} (\bsx_{a,t})-F_{t}(\bsx_t^*) \right] &\leq \sum\limits_{t=1}^{T-1} \mathbb{E} [\Omega_t] +2nL_X R + 2n H_T\nonumber \\
& \quad +  2 \alpha n R  D_T.
\end{flalign}

For the first term on the right hand of (\ref{proof lem key 6-a6}), we obtain
\begin{flalign}
 &\mathbb{E}[\Upsilon_{i,t}]=
  \mathbb{E}\left \{ \left\langle \nabla F_{t+1} (\boldsymbol{x}_{a,t}), \mathds{Q}(\bsx_{i,t})-\bsx_{i,t}  \right\rangle \right\}
=0,  \label{proof lem key 6-a7}\\
& {G_X}  \sum_{i=1}^n \mathbb{E}[\left\|\boldsymbol{e}_{i,t}\right\|^2]\leq nG_XR^2 \epsilon_{d,k_t}.  \label{proof lem key 6-a72}
\end{flalign}

Through recalling the definition of $\Omega_t$ and combining (\ref{Proof Theorem1- 1-a2}), (\ref{proof lem key 6-a6}), (\ref{proof lem key 6-a7}), (\ref{proof lem key 6-a72}), Lemmas \ref{consistency} and \ref{grad diffience}, we can readily obtain the condition (\ref{Eq Theorem1}).  The proof is complete.
\hfill$\square$

In Theorem 1, (\ref{Eq Theorem1}) shows that the upper bound of the regret depends on $D_i, i\in \{1,2,\ldots,6 \}, \alpha, T, \epsilon_{d,k_t}, H_T$ and $ D_T$, where $D_i$ are the scalars  consisting of the initial values, optimization problem parameters, and the network  parameters $\sigma, \Gamma$. It should be pointed out that the shorter the jointed connection period $Q$ of graph is, the tighter the regret bound will be due to the smaller parameter $\Gamma$ and $\frac{1}{1-\sigma}$ in coefficients $D_1, D_2, D_3 $ and $D_6$.
When the optimization problem and the network graph are determined, the coefficient $D_i$  will be a finite fixed constants and $D_T, H_T$ will have a fixed order over time $T$.
 It is not hard to note that the regret bound of Algorithm \ref{algorithm 1} is affected by $\alpha$ and $\epsilon_{d,k_t}$.
Thus, we have the following corollary.
\begin{corollary} \label{corollary 1}
Suppose that the conditions in Theorem \ref{theorem 1} and $H_T = o(T), D_T = o(T) $  hold. Given positive scalars $\kappa_1, \gamma<1,\kappa_2 \leq T^{\gamma}$ and $ \xi$,  the quantization resolution and the step size are chosen as $\epsilon_{d,k_t}=\frac{\kappa_1}{ t^{\xi}}, \alpha = \frac{\kappa_2}{T^{\gamma}}$, respectively. Then, we have that
\begin{flalign}\label{corollary 1 equation}
 &\mathbb{E}[\textbf{Regret}_d^j (T)] \leq \nonumber\\
&  \left \{
\begin{array}{l}
    \mathcal{O}\left(\max\left\{ T^{1-b},  T^{\gamma}(1+H_T) \right\} +D_T\right), \text{when} \  \gamma<\xi<1.\\
    \mathcal{O}\left( \max\left\{T^{1-\gamma}, T^{\gamma}\ln T, T^{\gamma}H_T \right\}   +D_T\right), \text{when} \  \xi=1.\\
    \mathcal{O}\left( \max\left\{T^{1-\gamma}, T^{\gamma}(1+H_T) \right\}   +D_T\right),  \text{when} \  1<\xi.\\
\end{array}
\right.
\end{flalign}
where $b:=\min\{\gamma,  \xi/2, \xi-\gamma\}$.
\end{corollary}
\noindent{\em Proof.} Substituting the conditions  of $\epsilon_{d,k_t}, \alpha$ in Corollary \ref{corollary 1} into (\ref{Eq Theorem1}), once can verify that  $\sum_{t=1}^T \sqrt{\epsilon_{d,k_t}}= \sum_{t=1}^T \frac{\sqrt{\kappa_1}}{t^{\xi/2}}=\sqrt{\kappa_1}+ \sqrt{\kappa_1} \int_1^T \frac{1}{t^{\xi/2}} dt \leq \mathcal{O}(T^{1-\xi/2})$ when $0<\xi<2$, $\sum_{t=1}^T \sqrt{\epsilon_{d,k_t}}\leq \mathcal{O}(\ln T)$ when $\xi=2$, and $\sum_{t=1}^T \sqrt{\epsilon_{d,k_t}}\leq \mathcal{O}(1)$ when $2< \xi$, respectively.
Similarly, $\sum_{t=1}^T \epsilon_{d,k_t}$ can be bounded as $ \mathcal{O}(T^{1-\xi})$ when $0<\xi<1$, $\mathcal{O}(\ln T)$ when $\xi=1$, and $ \mathcal{O}(1)$ when $1< \xi$, respectively.
Then, (\ref{corollary 1 equation}) is easily obtained based on different ranges of $\xi$.
The proof is complete.
\hfill$\square$

\begin{remark}
 Note that the setting of the decreasing quantization resolution in Corollary \ref{corollary 1} allows relatively coarse quantization in the early stage of algorithm execution. In particular,  when the parameter $\xi$ is chosen as a small value, the saving of communication resources will be significant but the regret bound will be poor, which implies that the setting of the parameter $\xi$ links the trade-off between them. It should be noted that when the total iteration time $T$ is large, the quantitative effect of information may be weakened to be close to the real value, especially in the later stage of algorithm operation. This change is actually reasonable
since the state variables must approach \emph{final} optima by continually obtaining the precise data as long as the algorithm runs.
\end{remark}
\begin{remark}
The result of Corollary \ref{corollary 1} matches the centralized result \cite{kalhan2021dynamic} and distributed results \cite{zhang2023dynamic, 10025380} while taking quantized communication into account. Compared with \cite{10025380}, we additionally consider  quantization communication and do not require the loss function to be bounded.
Moreover, the requirements $H_T = o(T)$ and $D_T = o(T)$ in Corollary \ref{corollary 1} imply that the cumulative variations of function value and gradient value grow slower than $T$ as $T$ increases.  This also means that the loss functions and gradient functions satisfy certain regularities over time, such as the variability of function parameters decreases over time.  According to Corollary \ref{corollary 1}, it is not hard to find that this requirement is reasonable and necessary for guaranteeing the sub-linearity of the considered dynamic regret.
In addition, if the bound of the prior knowledge $H_T$  can be known in advance, i.e., $H_T \leq \mathcal{O} (T^{\theta}), 0< \theta < 1 $, (\ref{corollary 1 equation}) can be improved to $\mathcal{O}( \sqrt{ T (1+H_T)} +D_T) $ by setting  $\gamma = 1/2 - \log_T \sqrt{1+T^{\theta}}$.
\end{remark}
\section{Simulation}
In this section, the following distributed online linear regression problem with a regularization term   is simulated to verify the proposed algorithm.
\begin{flalign}
&\min\limits_{\boldsymbol{x}\in \boldsymbol{X}} \sum_{t=1}^T \sum_{i=1}^n\left[\frac{1}{2}\left(\boldsymbol{p}_{i, t}^{\top} \bsx -q_{i, t}\right)^2+\rho\|\bsx\|_2^2\right]
\end{flalign}
where $\boldsymbol{X}: = \{ \bsx | \| \bsx\|_1 \leq 2 \}$, $\boldsymbol{p}_{i, t}  \in \mathbb{R}^d, q_{i, t} \in \mathbb{R}$ represents the feature and label information, and $\rho$ is a regular parameter. The feature vector $\boldsymbol{p}_{i, t}$ is generated randomly and  uniformly and its  element satisfies $[\boldsymbol{p}_{i, t}]_i \in [-5, 5]$. The label $q_{i, t}$ satisfies
 $
 q_{i, t}=\boldsymbol{p}_{i, t}^{\top} \bsx_0 + {\zeta_{i,t}}/(4 {t})
 $
 where $\zeta_{i,t}$ is generated randomly in the interval $[0, 1]$.
  In the following simulation, we set  the algorithm parameters $n=10$, $d=30$, $\rho= 5 \times 10^{-6}$, $\alpha=1 /(2T^{0.3})$ and take the probabilistic quantizer mentioned in Remark \ref{remark 1} as an example. To measure the performance of the algorithm, the global average dynamic regret $\frac{1}{n}\sum_{j=1}^n [\textbf{Regret}_d^j(T)/T] $ is defined.

\begin{figure}[ht]
 \centering
 \includegraphics[width=9cm]{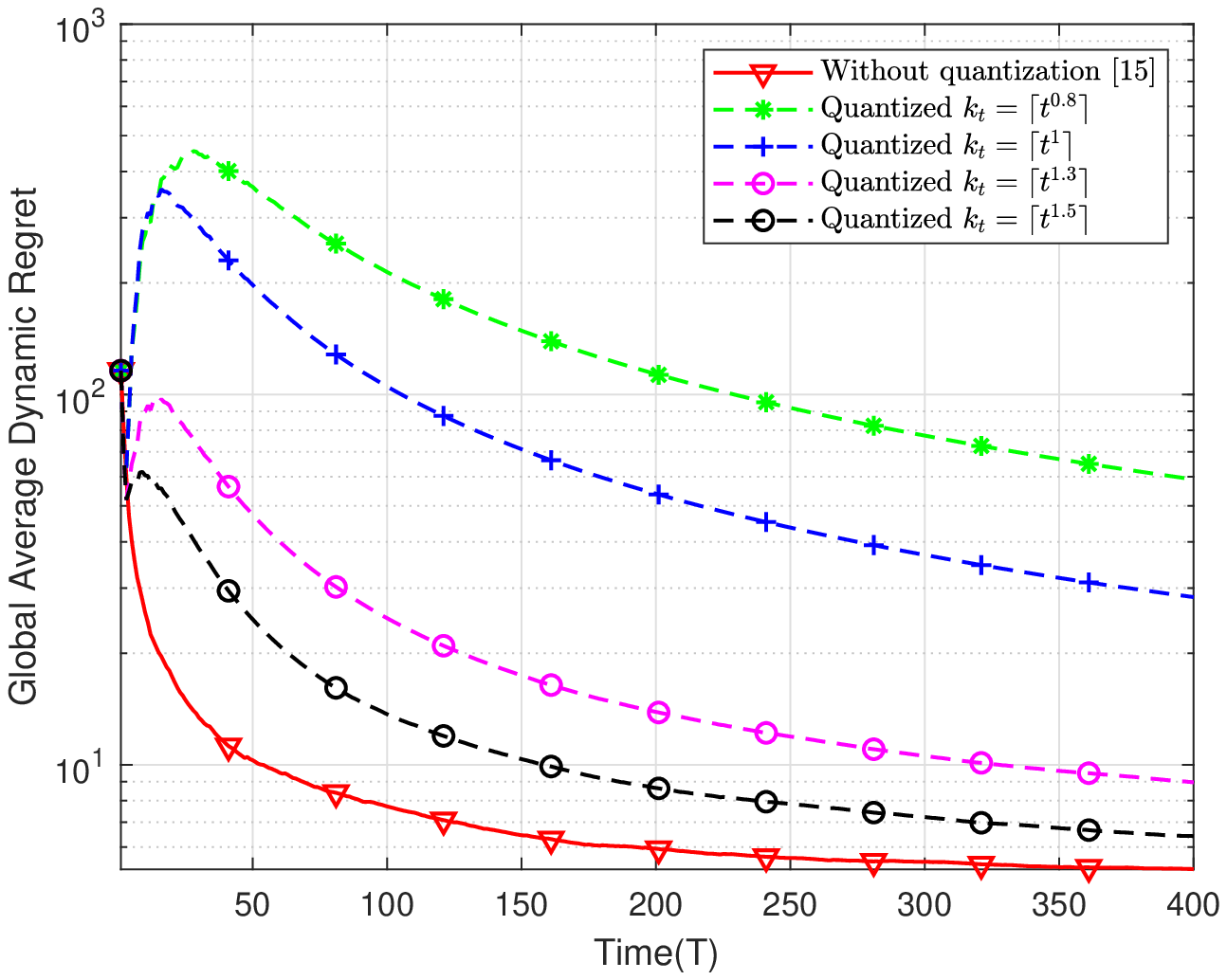}
  \caption{The comparison of Algorithm \ref{algorithm 1} under three quantization resolutions.}
  \label{f1}
   \includegraphics[width=9cm]{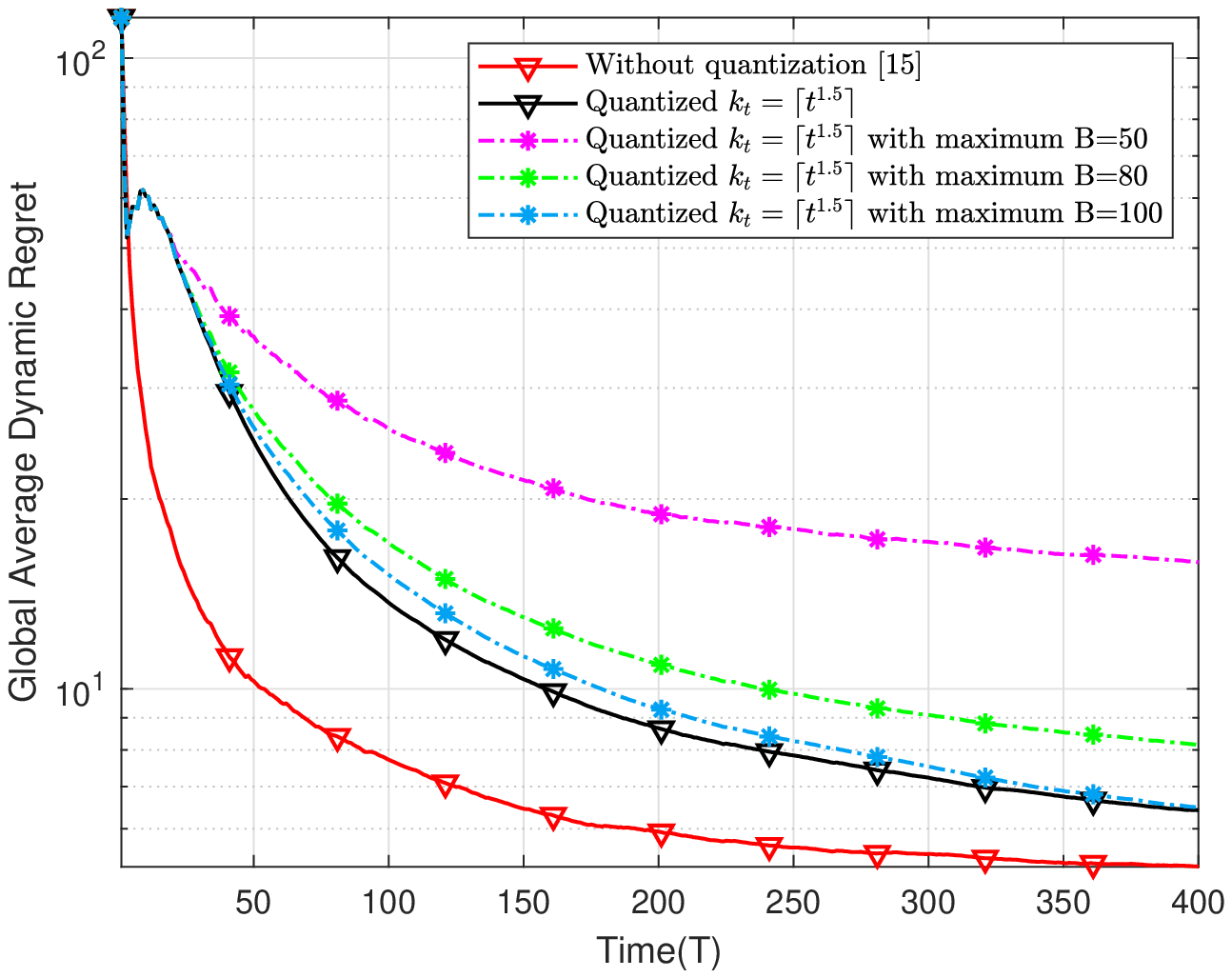}
  \caption{The effect of quantization level with an upper limit on convergence performance.}
  \label{f1_2}
  \end{figure}

To investigate the convergence of Algorithm \ref{algorithm 1} and the effect of quantization parameters on algorithm convergence, we compare the global average dynamic regrets of Algorithm \ref{algorithm 1} under different cases: no quantization \cite{zhang2023dynamic}, quantization levels $k_t= \lceil t^{0.8} \rceil, \lceil t^{1} \rceil , \lceil t^{1.3} \rceil$ and  $\lceil t^{1.5} \rceil$.
From Fig. \ref{f1},  Algorithm \ref{algorithm 1} is convergent and when $k_t$ with a larger increasing tendency is selected, the related convergence performance  is better. Note that in the early stage of the algorithm, the performance fluctuation caused by relatively coarse quantization resolution can be tolerated and this error can be weakened with the iteration time.
Further, we analyze the effect of the quantization level with the maximum number $B \geq k_t$ on the convergence performance of the designed algorithm. Taking quantization level $k_t=\lceil t^{1.5} \rceil$ as an example and considering the case of $B=50, 80, 100,$ the comparison results are shown in Fig. \ref{f1_2}. It can be seen that when $B=100$, its convergence curve is close to that without the maximum number, while $B=50$, the convergence performance is poor. It should be noted that the design of quantization level with an appropriate parameter $B$ can better save communication resources than that without the limited parameter $B$, but it always has a quantization error $\boldsymbol{e}_{i,t}$ because the quantized data cannot approach the real-value data.

\begin{figure}[ht]
 \includegraphics[width=9cm]{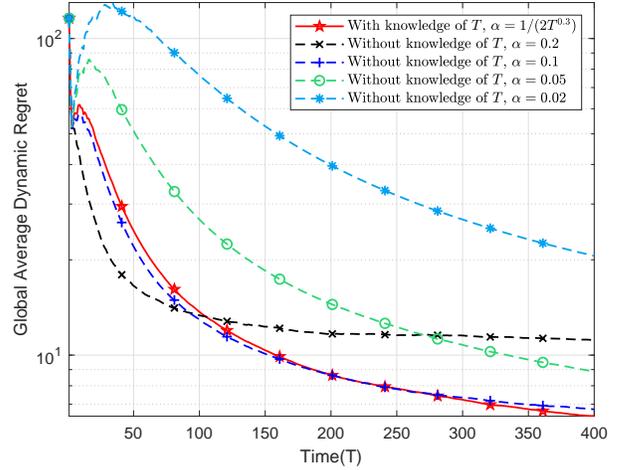}
  \caption{The performance comparison under the step sizes with and without the knowledge of $T$.}
  \label{f1_3}
  \end{figure}

\begin{figure}[ht]
   \includegraphics[width=9cm]{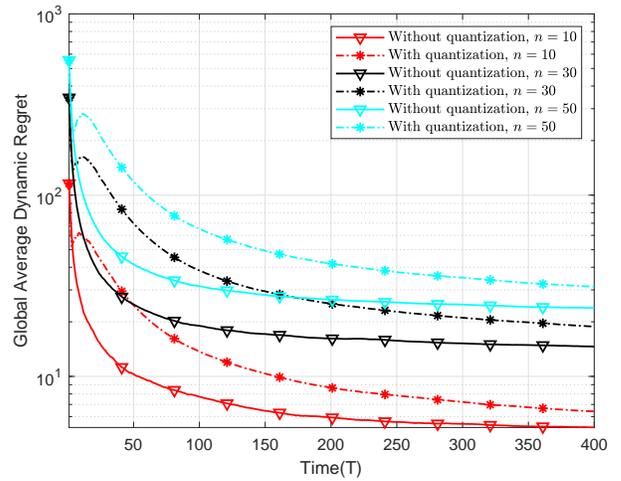}
  \caption{The effect of the number of agents on the convergence performance.}
  \label{f2}
\end{figure}

Next, we carry out a comparative study for the convergence performance of Algorithm \ref{algorithm 1} under the step size design taking into account unknown total iteration time $T$ as well as the quantization level $k_t=\lceil t^{1.5} \rceil$. Without loss of generality, setting the step size as $\alpha= 1/(2 T^{0.3}),  0.2,0.1,0.05,0.02$, respectively, the comparison results of the convergence performance are revealed in Fig. \ref{f1_3}. Among the settings of step sizes, the dynamic regret under the step size with the knowledge of $T$ has a significantly better convergence effect,  while that under the step size without the knowledge of $T$ has a large fluctuation of the convergence performance for different settings. Although the latter does not require prior knowledge of $T$, it always has a performance gap in a theoretical sense according to Theorem \ref{theorem 1}. In addition, as the horizon $T$ varies,  this step size without the horizon $T$ may cause the original convergence performance to be unmaintainable due to its \emph{invariant setting}.
Finally, the effect of the number of agents on the convergence performance is studied under $k_t= \lceil t^{1.5} \rceil$. Through setting $n=10, 30, 50$, the comparison of the global average dynamic regret is shown in Fig. \ref{f2}, which verifies the theoretical results in Theorem \ref{theorem 1} that the smaller the value of $n$ is, the better the convergence performance of Algorithm \ref{algorithm 1} is.

\section{CONCLUSIONS}
For the distributed online constrained optimization problem under quantized communication, this paper has developed a quantized
distributed online projection-free optimization algorithm. The use of random quantizers and linear oracle in the proposed algorithm has ensured the effective saving of communication resources and computational costs.  For different settings of quantization resolution  $\epsilon_{d,k_t}$,
the related dynamic regret bound  has been established, in which the optimal bound $\mathcal{O}( \sqrt{ T (1+H_T)} +D_T) $ can be achieved when the knowledge $H_T$ is known and $\xi>1$. In addition, we have revealed the trade-off between the convergence performance and the quantization effect. Finally, a simulation example has been investigated to verify the theoretical results. A promising direction in the future is to investigate the nonconvex loss function case that will be more general yet more challenging.

%

\appendix


\subsection{Proof of lemma \ref{consistency}}
According to Algorithm \ref{algorithm 1}, we get
\begin{flalign} \label{proof_consistency 1}
\hat{\bsx}_{i,t}
&=\sum\limits_{j=1}^n [W_t]_{ij}  \mathds{Q}_t(\bsx_{j,t})\nonumber \\
&=  \sum\limits_{j=1}^n [W_t]_{ij}  \left[ \hat{\bsx}_{j,t-1} +\alpha (\boldsymbol{v}_{j,t-1}-\hat{\bsx}_{j,t-1} )\right] \nonumber \\
&\quad + \sum\limits_{j=1}^n [W_t]_{ij} \bse_{j,t}\nonumber \\
%
&= \sum\limits_{j=1}^n [\Phi(t,2)]_{ij}  \hat{\bsx}_{j,1} + \sum\limits_{l=2}^{t} \sum\limits_{j=1}^n [\Phi(t,l)]_{ij} \boldsymbol{e}_{j,l} \nonumber \\
&\quad +\alpha \sum\limits_{l=1}^{t-1} \sum\limits_{j=1}^n [\Phi(t,l+1)]_{ij} (\boldsymbol{v}_{j,l}-\hat{\bsx}_{j,l} ) \nonumber \\
&= \sum\limits_{j=1}^n [\Phi(t,1)]_{ij}  {\bsx}_{j,1} + \sum\limits_{l=1}^{t} \sum\limits_{j=1}^n [\Phi(t,l)]_{ij} \boldsymbol{e}_{j,l} \nonumber \\
&\quad +\alpha \sum\limits_{l=1}^{t-1} \sum\limits_{j=1}^n [\Phi(t,l+1)]_{ij} (\boldsymbol{v}_{j,l}-\hat{\bsx}_{j,l} ).
\end{flalign}

According to Algorithm \ref{algorithm 1}, the term ${\boldsymbol{x}}_{a,t}$ can be further simplified as follows:
\begin{flalign} \label{proof_consistency 2}
&{\boldsymbol{x}}_{a,t}\nonumber \\
&=\frac{1}{n} \sumn \sum_{j=1}^n [W_{t-1}]_{ij} \mathds{Q}_{t-1}(\boldsymbol{x}_{j,t-1})+ \frac{\alpha}{n}\sum_{i=1}^n \left({\boldsymbol{v}}_{i,t-1}-\hat{\boldsymbol{x}}_{i,t-1}\right)\nonumber \\
&=\frac{1}{n} \sumn \mathds{Q}_{t-1}(\boldsymbol{x}_{i,t-1})+ \frac{\alpha}{n}\sum_{i=1}^n \left({\boldsymbol{v}}_{i,t-1}-\hat{\boldsymbol{x}}_{i,t-1}\right)\nonumber \\
&={\boldsymbol{x}}_{a,t-1}+ \frac{1}{n} \sumn \bse_{i,t-1} + \frac{\alpha}{n}\sum_{i=1}^n \left({\boldsymbol{v}}_{i,t-1}-\hat{\boldsymbol{x}}_{j,t-1}\right)\nonumber \\
%
&= \frac{1}{n} \sum\limits_{j=1}^n {\boldsymbol{x}}_{j,1}+ \frac{1}{n}\sum_{l=1}^{t-1} \sumn \bse_{i,l} +\frac{\alpha}{n} \sum\limits_{l=1}^{t-1}  \sum\limits_{j=1}^n ({\boldsymbol{v}}_{j,l}-\hat{\boldsymbol{x}}_{j,l})
\end{flalign}
where the second equality combines the double stochasticity of weight matrix
$W_{t-1}$.

 Combining   (\ref{proof_consistency 1}) and (\ref{proof_consistency 2}), for $t \geq 2$, we achieve
\begin{flalign} \label{proof_consistency 3}
& \| \hat{\bsx}_{i,t}- \bsx_{a,t}\|\nonumber \\
%
%
&\leq  \sum\limits_{j=1}^n \left| [\Phi(t,1)]_{ij}  - \frac{1}{n} \right| \| {\boldsymbol{x}}_{j,1} \| \nonumber \\
&\quad +\sum\limits_{l=1}^{t-1} \sum\limits_{j=1}^n  \left|[\Phi(t,l)]_{ij} - \frac{1}{n}\right| \|\bse_{j,l} \| +\sum_{j=1}^n [W_t]_{ij} \|\bse_{j,t}\|\nonumber \\
&\quad +\alpha \sum\limits_{l=1}^{t-1} \sum\limits_{j=1}^n  \left|[\Phi(t,l+1)]_{ij} - \frac{1}{n}\right|   \left\| {\boldsymbol{v}}_{j,l}-\hat{\boldsymbol{x}}_{j,l}\right\|\nonumber \\
&\leq \Gamma \sigma^{t-1}\sum\limits_{j=1}^n  \| {\boldsymbol{x}}_{j,1} \|+2\alpha nR\Gamma \sum\limits_{l=1}^{t-1}  \sigma^{t-l-1} \nonumber \\
&\quad +\sum_{j=1}^n [W_t]_{ij} \|\bse_{j,t}\| +\Gamma \sum\limits_{l=1}^{t-1} \sumn  \sigma^{t-l} \|\bse_{i,l} \|
\end{flalign}
where the second inequality follows the fact $\boldsymbol{v}_{j,l}, \hat{\boldsymbol{x}}_{j,l} \in  \boldsymbol{X}$ and the property \cite{nedic2008distributed} : $ \left|[\Phi(t, s)]_{i j}-\frac{1}{n}\right| \leq \Gamma \sigma^{(t-s)} $, $ \forall i, j \in \mathcal{V}$, $\sigma=(1-\zeta / 4 n^{2})^{1 /  Q},  \Gamma=(1-\zeta / 4 n^{2})^{(1-2Q) /  Q}$.

Summing from $i=1$ to $n$ and $t=1$ to $T$ on both sides of  (\ref{proof_consistency 3}), we get
\begin{flalign} \label{proof_consistency 4}
 &\sumT \sumn \mathbb{E} [\| \hat{\bsx}_{i,t}- \bsx_{a,t}\|] \nonumber \\
 &\leq \sumn \mathbb{E} [\| \hat{\bsx}_{i,1}- \bsx_{a,1}\|]  +\sum_{t=2}^T \sum_{i=1}^n  \mathbb{E} [\|\bse_{i,t}\|] +\alpha T\frac{2n^2R\Gamma}{1-\sigma} \nonumber \\
 &\quad +n\Gamma \sum_{t=2}^T \sum\limits_{l=1}^{t-1} \sumn  \sigma^{t-l} \mathbb{E} [\|\bse_{i,l} \|] + n\Gamma \sum_{t=2}^T \sigma^{t-1}\sum\limits_{j=1}^n  \| {\boldsymbol{x}}_{j,1} \| \nonumber \\
  &\leq    \frac{n\Gamma}{1-\sigma} \sum\limits_{j=1}^n  \| {\boldsymbol{x}}_{j,1} \|+\left(1+\frac{n \Gamma \sigma}{1-\sigma}\right)\sum_{t=1}^T \sum_{i=1}^n  \mathbb{E} [\|\bse_{i,t}\|]  \nonumber \\
  &\quad +\alpha T\frac{2n^2R\Gamma}{1-\sigma}
\end{flalign}
where the second inequality follows the fact
 \begin{flalign}
& \sum_{i=1}^n \| \hat{\bsx}_{i,1}- \bsx_{a,1}\| \nonumber \\
&\leq \sum_{i=1}^n \left\| \sum_{j=1}^n [W_1]_{ij} ({\bsx}_{j,1}+\bse_{j,1})- \frac{1}{n} \sum_{j=1}^n \bsx_{j,1}\right\| \nonumber \\
&\leq \sum_{i=1}^n \sum_{j=1}^n| [W_1]_{ij}- \frac{1}{n}| \|  \bsx_{j,1}\| +\sum_{i=1}^n \| \bse_{i,1}\|\nonumber \\
&\leq n\Gamma  \sum_{j=1}^n  \|  \bsx_{j,1}\| + \sum_{i=1}^n \| \bse_{i,1}\|.
\end{flalign}
The proof is complete.
\hfill$\square$
\subsection{Proof of lemma \ref{grad diffience}}
According to Assumption \ref{assump: lips grad}, we establish that
\begin{flalign} \label{proof_grad diff 0}
&  \left\| \frac{1}{n}  {\nabla} F_{t} ({\boldsymbol{x}}_{a,t})- \widehat{\bss} _{i,t}  \right\| \nonumber \\
& =  \left\|\frac{1}{n} \sum_{j=1}^n   {\nabla} f_{j,t} (\hat{\boldsymbol{x}}_{j,t}) - \widehat{\bss} _{i,t} \right. \nonumber\\
&\quad \left.- \frac{1}{n} \sum_{j=1}^n   [{\nabla} f_{j,t} (\hat{\boldsymbol{x}}_{j,t}) -{\nabla} f_{j,t} ({\boldsymbol{x}}_{a,t})] \right\| \nonumber \\
& \leq  \left\|\frac{1}{n} \sum_{j=1}^n   {\nabla} f_{j,t} (\hat{\boldsymbol{x}}_{j,t}) - \widehat{\bss} _{i,t}\right\| +\frac{G_X}{n} \sum_{j=1}^n  \left\| \hat{\boldsymbol{x}}_{j,t}- {\boldsymbol{x}}_{a,t} \right\|.
\end{flalign}
For the first term on the right hand side of (\ref{proof_grad diff 0}), from
 Algorithm \ref{algorithm 1},  for any  $t\geq2$ it can be verified that
\begin{flalign} \label{proof_grad diff 1}
\widehat{\bss} _{i,t} &=\sum\limits_{j=1}^n [W_t]_{ij}  \overline{\nabla} f_{j,t}\nonumber \\
&= \sum\limits_{j=1}^n [W_t]_{ij} \widehat{\bss} _{j,t-1} + \sum\limits_{j=1}^n [W_t]_{ij} \boldsymbol{\nabla}_{i,t}^Q\nonumber \\
&= \sum\limits_{j=1}^n [\Phi(t,2)]_{ij}\widehat{\bss} _{j,1} + \sum\limits_{l=2}^t \sum\limits_{j=1}^n [\Phi(t,l)]_{ij} \boldsymbol{\nabla}_{i,l}^Q \nonumber \\
&= \sum\limits_{j=1}^n [\Phi(t,1)]_{ij} \mathds{Q}_1 [{\nabla} f_{j,1} (\hat{\boldsymbol{x}}_{j,1})] + \sum\limits_{l=2}^t \sum\limits_{j=1}^n [\Phi(t,l)]_{ij} \boldsymbol{\nabla}_{i,l}^Q.\nonumber \\
\end{flalign}

Note that $\overline{\nabla} f_{i,1} = \mathds{Q}_1[\nabla f_{i,1} (\hat{\boldsymbol{x}}_{i,1})]$ from Algorithm \ref{algorithm 1}. Hence,   the equality $\sum_{i=1}^n \overline{\nabla} f_{i,t} = \sum_{i=1}^n \mathds{Q}_t[\nabla f_{i,t} (\hat{\boldsymbol{x}}_{i,t})]$  holds when $t=1$. Now  we assume that $\sum_{i=1}^n \overline{\nabla} f_{i,t-1} = \sum_{i=1}^n \mathds{Q}_{t-1}[\nabla f_{i,t-1} (\hat{\boldsymbol{x}}_{i,t-1})]$ holds at time $t-1$, and we intend to show the same conclusion at time $t$. Actually,
%
%
\begin{flalign}
 \sumn \overline{\nabla} f_{i,t}
& = \sumn \widehat{\bss} _{i,t-1} + \sumn  \mathds{Q}_t[\nabla f_{i,t}  (\hat{\boldsymbol{x}}_{i,t})] \nonumber \\
&\quad - \sumn  \mathds{Q}_{t-1}[\nabla f_{i,t-1}  (\hat{\boldsymbol{x}}_{i,t-1})] \nonumber \\
&= \sumn \sum_{j=1}^n [W_{t-1}]_{ij} \overline{\nabla} f_{j,t-1} + \sumn  \mathds{Q}_t[\nabla f_{i,t}  (\hat{\boldsymbol{x}}_{i,t})]\nonumber \\
&\quad  - \sum_{i=1}^n \overline{\nabla} f_{i,t-1}  \nonumber \\
&= \sumn  \mathds{Q}_t[\nabla f_{i,t}  (\hat{\boldsymbol{x}}_{i,t})]
\end{flalign}
where the last equality follows from the double stochasticity of $W_{t}$.
With this condition, we obtain for any $t\geq 2$ that
\begin{flalign} \label{proof_grad diff 2}
&\frac{1}{n} \sumn   {\nabla} f_{i,t} (\hat{\boldsymbol{x}}_{i,t}) +\frac{1}{n} \sumn \boldsymbol{\theta}_{i,t}\nonumber \\
&=  \frac{1}{n} \sumn   \overline{\nabla} f_{i,t}\nonumber \\
&= \frac{1}{n} \sumn \sum\limits_{j=1}^n [W_{t-1}]_{ij}  \overline{\nabla} f_{j,t-1} +  \frac{1}{n} \sumn  \boldsymbol{\nabla}_{i,t}^Q   \nonumber \\
&= \frac{1}{n} \sumn   \overline{\nabla} f_{i,t-1} +  \frac{1}{n} \sumn  \boldsymbol{\nabla}_{i,t}^Q    \nonumber \\
&= \frac{1}{n} \sumn   \overline{\nabla} f_{i,1} + \frac{1}{n} \sum\limits_{l=2}^t\sumn   \boldsymbol{\nabla}_{i,l}^Q   \nonumber \\
&= \frac{1}{n} \sumn   \mathds{Q}_1[{\nabla} f_{i,1} (\hat{\boldsymbol{x}}_{i,1})] + \frac{1}{n}\sum\limits_{l=2}^t \sumn  \boldsymbol{\nabla}_{i,l}^Q.
\end{flalign}

Similar to (\ref{proof_consistency 3}), combining the fact $\widehat{\bss} _{i,1}=\sum_{j=1}^n [W_1]_{ij}  \mathds{Q}_1 [{\nabla} f_{j,1} (\hat{\boldsymbol{x}}_{j,1})]$, it follows from (\ref{proof_grad diff 1}) and (\ref{proof_grad diff 2}) that
\begin{flalign}
&\|  \widehat{\bss} _{i,t} - \frac{1}{n} \sum_{j=1}^n  {\nabla} f_{j,t} (\hat{\boldsymbol{x}}_{j,t})\|  \nonumber \\
&\leq\sum_{j=1}^n \Gamma \sigma^{t-1}  \| \mathds{Q}_1[  {\nabla} f_{j,1} (\hat{\boldsymbol{x}}_{j,1})] \| +\sum_{l=2}^t \sum_{j=1}^n\Gamma \sigma^{t-l} \|  \boldsymbol{\nabla}_{j,l}^Q \| \nonumber \\
&\quad+ \frac{1}{n}\sum_{i=1}^n \|\boldsymbol{\theta}_{i,t}\|.
\end{flalign}
This implies that
\begin{flalign} \label{proof_grad diff 4}
 &\sumT \sumn \mathbb{E}\left[\left\|  \widehat{\bss} _{i,t} - \frac{1}{n} \sum_{j=1}^n   {\nabla} f_{j,t} (\hat{\boldsymbol{x}}_{j,t})\right\|\right]\nonumber \\
 &\leq\sumn \mathbb{E}\left[\left\|  \widehat{\bss} _{i,1} -\frac{1}{n} \sum_{j=1}^n    {\nabla} f_{j,1} (\hat{\boldsymbol{x}}_{j,1})\right\| \right] \nonumber \\
 &+\sum\limits_{t=2}^T \sum\limits_{j=1}^n n \Gamma \sigma^{t-1}  \mathbb{E}\left\{\left\| \mathds{Q}_1[  {\nabla} f_{j,1} (\hat{\boldsymbol{x}}_{j,1}) ]\right\| \right\}\nonumber \\
  &\quad +\sum\limits_{t=2}^T \sum\limits_{l=2}^t \sum\limits_{j=1}^n n \Gamma \sigma^{t-l} \mathbb{E}\left[\left\|  \boldsymbol{\nabla}_{j,l}^Q  \right\| \right] +\frac{1}{n} \sum_{t=2}^T\sum_{i=1}^n \mathbb{E}\left[\|\boldsymbol{\theta}_{i,t}\| \right]\nonumber \\
 &\leq\sumn \sum\limits_{j=1}^n \left | [W_1]_{ij} -  \frac{1}{n} \right| \left\|  {\nabla} f_{j,1} (\hat{\boldsymbol{x}}_{j,1}) \right\| +\sumn \mathbb{E}\left[\|\boldsymbol{\theta}_{i,1}\|\right] \nonumber \\
 &+\frac{\sigma n \Gamma}{1-\sigma} \sum\limits_{j=1}^n \left\|   {\nabla} f_{j,1} (\hat{\boldsymbol{x}}_{j,1}) \right\| + \frac{\sigma n \Gamma}{1-\sigma} \sumn \mathbb{E}\left[\|\boldsymbol{\theta}_{i,1}\|\right]\nonumber \\
 & + n\Gamma \left(\sum\limits_{l=0}^{T-2}
 \sigma^l \right)\left(\sum\limits_{t=2}^{T} \sumn \mathbb{E}\left[\left\|  \boldsymbol{\nabla}_{i,t}^Q  \right\| \right]\right)+ \sum_{t=2}^T\sum_{i=1}^n \mathbb{E}\left[\|\boldsymbol{\theta}_{i,t}\| \right] \nonumber \\
 &\leq \frac{n \Gamma}{1-\sigma} \sum\limits_{j=1}^n \left\|   {\nabla} f_{j,1} (\hat{\boldsymbol{x}}_{j,1}) \right\|+\frac{ n \Gamma}{1-\sigma} \sum\limits_{t=2}^{T} \sumn \mathbb{E}\left[\left\|   \boldsymbol{\nabla}_{i,t}^Q  \right\|\right]\nonumber \\
 &\quad + \frac{\sigma n \Gamma}{1-\sigma} \sumn \mathbb{E}\left[\|\boldsymbol{\theta}_{i,1}\|\right] + \sum_{t=1}^T\sum_{i=1}^n \mathbb{E}\left[\|\boldsymbol{\theta}_{i,t}\| \right].
\end{flalign}

For the second term on the right hand side of (\ref{proof_grad diff 4}), by recalling the notion $ \boldsymbol{\nabla}_{i,t}^Q$ defined in (\ref{average-delta}),  it can be yielded that
\begin{flalign} \label{proof Sigma Tracking 1}
 &\sumn \left\|   \boldsymbol{\nabla}_{i,t}^Q \right\| -\sumn ( \|\boldsymbol{\theta}_{i,t}\| + \|\boldsymbol{\theta}_{i,t-1}\| )  \nonumber \\
 & \leq \sumn \left\| \nabla f_{i,t}  (\hat{\bsx}_{i,t}) - \nabla f_{i,t-1}  (\hat{\boldsymbol{x}}_{i,t})\right\| \nonumber \\
 &\quad +\sumn \left\| \nabla f_{i,t-1}  (\hat{\boldsymbol{x}}_{i,t}) - \nabla f_{i,t-1}  (\hat{\boldsymbol{x}}_{i,t-1})\right\| \nonumber \\
 &\leq n g_{t-1}^{sup} +\sumn \left(  G_{X} \| \hat{\boldsymbol{x}}_{i,t} -\hat{\boldsymbol{x}}_{i,t-1} \|  \right)\nonumber \\
 &\leq n g_{t-1}^{sup} +G_{X}\sumn (  \| \hat{\boldsymbol{x}}_{i,t} -{\boldsymbol{x}}_{a,t-1} \|   +\| \hat{\boldsymbol{x}}_{i,t-1}-{\boldsymbol{x}}_{a,t-1} \| )\nonumber \\
 &\leq n g_{t-1}^{sup} + 2G_{X}\sumn \| \hat{\boldsymbol{x}}_{i,t-1} -{\boldsymbol{x}}_{a,t-1} \|+G_X\sum_{i=1}^n \| \bse_{i,t}\|   \nonumber \\
 & \quad +2n R G_{X} \alpha
\end{flalign}
where the last inequality is obtained based on the fact:
\begin{flalign}
&\sum_{i=1}^n \| \hat{\boldsymbol{x}}_{i,t} -{\boldsymbol{x}}_{a,t-1} \| \nonumber \\
&\leq \sum_{i=1}^n \sum_{j=1}^n [W_t]_{ij}\| \mathds{Q}({\boldsymbol{x}}_{j,t}) -{\boldsymbol{x}}_{a,t-1} \| \nonumber \\
&\leq \sum_{i=1}^n \| \bse_{i,t}+\hat{\boldsymbol{x}}_{i,t-1} -{\boldsymbol{x}}_{a,t-1} +\alpha (\boldsymbol{v}_{i,t-1}-\hat{\boldsymbol{x}}_{i,t-1})\| \nonumber \\
&\leq \sum_{i=1}^n \| \hat{\boldsymbol{x}}_{i,t-1} -{\boldsymbol{x}}_{a,t-1} \| +\sum_{i=1}^n \| \bse_{i,t}\|+ 2\alpha n R.
\end{flalign}

Substituting the above inequalities into   (\ref{proof_grad diff 4}) and combining (\ref{proof_grad diff 0}) and fact $\mathbb{E} [\| \boldsymbol{\theta}_{i,t}\|]  \leq \sqrt{\mathbb{E} [\| \boldsymbol{\theta}_{i,t}\|^2]} \leq \sqrt{\epsilon_{d,k_t} \| \nabla f_{i,t}  (\hat{\bsx}_{i,t})\|^2  } \leq L_X \sqrt{\epsilon_{d,k_t}}$,  we can readily obtain (\ref{condition-lem4}). The proof is complete.
\hfill$\square$

\bibliographystyle{IEEEtran}
\bibliography{V3Arxiv_Quantized_DOPF}

\begin{thebibliography}{10}
\providecommand{\url}[1]{#1}
\csname url@samestyle\endcsname
\providecommand{\newblock}{\relax}
\providecommand{\bibinfo}[2]{#2}
\providecommand{\BIBentrySTDinterwordspacing}{\spaceskip=0pt\relax}
\providecommand{\BIBentryALTinterwordstretchfactor}{4}
\providecommand{\BIBentryALTinterwordspacing}{\spaceskip=\fontdimen2\font plus
\BIBentryALTinterwordstretchfactor\fontdimen3\font minus
  \fontdimen4\font\relax}
\providecommand{\BIBforeignlanguage}[2]{{%
\expandafter\ifx\csname l@#1\endcsname\relax
\typeout{** WARNING: IEEEtran.bst: No hyphenation pattern has been}%
\typeout{** loaded for the language `#1'. Using the pattern for}%
\typeout{** the default language instead.}%
\else
\language=\csname l@#1\endcsname
\fi
#2}}
\providecommand{\BIBdecl}{\relax}
\BIBdecl

\bibitem{yi2020distributed}
X.~Yi, X.~Li, T.~Yang, L.~Xie, T.~Chai, and K.~H. Johansson, ``Distributed
  bandit online convex optimization with time-varying coupled inequality
  constraints,'' \emph{IEEE Transactions on Automatic Control}, vol.~66,
  no.~10, pp. 4620--4635, 2021.

\bibitem{shahrampour2017distributed}
S.~Shahrampour and A.~Jadbabaie, ``Distributed online optimization in dynamic
  environments using mirror descent,'' \emph{IEEE Transactions on Automatic
  Control}, vol.~63, no.~3, pp. 714--725, 2017.

\bibitem{yuan2022distributedauto}
D.~Yuan, B.~Zhang, D.~W. Ho, W.~X. Zheng, and S.~Xu, ``Distributed online
  bandit optimization under random quantization,'' \emph{Automatica}, vol. 146,
  p. 110590, 2022.

\bibitem{yang2019survey}
T.~Yang, X.~Yi, J.~Wu, Y.~Yuan, D.~Wu, Z.~Meng, Y.~Hong, H.~Wang, Z.~Lin, and
  K.~H. Johansson, ``A survey of distributed optimization,'' \emph{Annual
  Reviews in Control}, vol.~47, pp. 278--305, 2019.

\bibitem{li2022survey}
X.~Li, L.~Xie, and N.~Li, ``A survey of decentralized online learning,''
  \emph{arXiv preprint}, vol. arXiv:2205.00473, 2022.

\bibitem{8715380}
C.~Liu, H.~Li, Y.~Shi, and D.~Xu, ``Distributed event-triggered gradient method
  for constrained convex minimization,'' \emph{IEEE Transactions on Automatic
  Control}, vol.~65, no.~2, pp. 778--785, 2020.

\bibitem{sundhar2010distributed}
S.~Sundhar~Ram, A.~Nedi{\'c}, and V.~V. Veeravalli, ``Distributed stochastic
  subgradient projection algorithms for convex optimization,'' \emph{Journal of
  Optimization Theory and Applications}, vol. 147, no.~3, pp. 516--545, 2010.

\bibitem{6311406}
F.~Yan, S.~Sundaram, S.~Vishwanathan, and Y.~Qi, ``Distributed autonomous
  online learning: Regrets and intrinsic privacy-preserving properties,''
  \emph{IEEE Transactions on Knowledge and Data Engineering}, vol.~25, no.~11,
  pp. 2483--2493, 2013.

\bibitem{pmlr-v70-zhang17g}
W.~Zhang, P.~Zhao, W.~Zhu, S.~C.~H. Hoi, and T.~Zhang, ``Projection-free
  distributed online learning in networks,'' in \emph{Proceedings of the 34th
  International Conference on Machine Learning}, 2017, pp. 4054--4062.

\bibitem{hazan2012projection}
E.~Hazan and S.~Kale, ``Projection-free online learning,'' in \emph{Proceedings
  of the 29th International Coference on International Conference on Machine
  Learning}, 2012, pp. 1843--1850.

\bibitem{7883821}
H.-T. Wai, J.~Lafond, A.~Scaglione, and E.~Moulines, ``Decentralized
  {Frank-Wolfe} algorithm for convex and nonconvex problems,'' \emph{IEEE
  Transactions on Automatic Control}, vol.~62, no.~11, pp. 5522--5537, 2017.

\bibitem{pmlr-v119-wan20b}
Y.~Wan, W.-W. Tu, and L.~Zhang, ``Projection-free distributed online convex
  optimization with $o(\sqrt{T})$ communication complexity,'' in
  \emph{Proceedings of the 37th International Conference on Machine Learning},
  2020, pp. 9818--9828.

\bibitem{wan2021projection}
Y.~Wan, G.~Wang, and L.~Zhang, ``Projection-free distributed online learning
  with strongly convex losses,'' \emph{arXiv preprint}, vol. arXiv:2103.11102,
  2021.

\bibitem{thuang2022stochastic}
N.~K. Thang, A.~Srivastav, D.~Trystram, and P.~Youssef, ``A stochastic
  conditional gradient algorithm for decentralized online convex
  optimization,'' \emph{Journal of Parallel and Distributed Computing}, vol.
  169, pp. 334--351, 2022.

\bibitem{zhang2023dynamic}
W.~Zhang, Y.~Shi, B.~Zhang, and D.~Yuan, ``Dynamic regret of distributed online
  frank-wolfe convex optimization,'' \emph{arXiv preprint arXiv:2302.00663},
  2023.

\bibitem{10025380}
K.~Lu and L.~Wang, ``Online distributed optimization with nonconvex objective
  functions via dynamic regrets,'' \emph{IEEE Transactions on Automatic
  Control}, 2023, doi: 10.1109/TAC.2023.3239432.

\bibitem{cao2023communication}
X.~Cao, T.~Ba{\c{s}}ar, S.~Diggavi, Y.~C. Eldar, K.~B. Letaief, H.~V. Poor, and
  J.~Zhang, ``Communication-efficient distributed learning: An overview,''
  \emph{IEEE Journal on Selected Areas in Communications}, 2023, doi:
  10.1109/JSAC.2023.3242710.

\bibitem{xiong2022distributed}
M.~Xiong, B.~Zhang, D.~Yuan, and S.~Xu, ``Distributed quantized mirror descent
  for strongly convex optimization over time-varying directed graph,''
  \emph{Science China Information Sciences}, vol.~65, no.~10, p. 202202, 2022.

\bibitem{yi2014quantized}
P.~Yi and Y.~Hong, ``Quantized subgradient algorithm and data-rate analysis for
  distributed optimization,'' \emph{IEEE Transactions on Control of Network
  Systems}, vol.~1, no.~4, pp. 380--392, 2014.

\bibitem{huang2016quantized}
C.~Huang, H.~Li, D.~Xia, and L.~Xiao, ``Quantized subgradient algorithm with
  limited bandwidth communications for solving distributed optimization over
  general directed multi-agent networks,'' \emph{Neurocomputing}, vol. 185, pp.
  153--162, 2016.

\bibitem{pu2016quantization}
Y.~Pu, M.~N. Zeilinger, and C.~N. Jones, ``Quantization design for distributed
  optimization,'' \emph{IEEE Transactions on Automatic Control}, vol.~62,
  no.~5, pp. 2107--2120, 2016.

\bibitem{9224135}
T.~T. Doan, S.~T. Maguluri, and J.~Romberg, ``Convergence rates of distributed
  gradient methods under random quantization: A stochastic approximation
  approach,'' \emph{IEEE Transactions on Automatic Control}, vol.~66, no.~10,
  pp. 4469--4484, 2021.

\bibitem{7891027}
H.~Li, C.~Huang, G.~Chen, X.~Liao, and T.~Huang, ``Distributed consensus
  optimization in multiagent networks with time-varying directed topologies and
  quantized communication,'' \emph{IEEE Transactions on Cybernetics}, vol.~47,
  no.~8, pp. 2044--2057, 2017.

\bibitem{9157925}
T.~T. Doan, S.~T. Maguluri, and J.~Romberg, ``Fast convergence rates of
  distributed subgradient methods with adaptive quantization,'' \emph{IEEE
  Transactions on Automatic Control}, vol.~66, no.~5, pp. 2191--2205, 2021.

\bibitem{9226092}
S.~Magn\'{u}sson, H.~Shokri-Ghadikolaei, and N.~Li, ``On maintaining linear
  convergence of distributed learning and optimization under limited
  communication,'' \emph{IEEE Transactions on Signal Processing}, vol.~68, pp.
  6101--6116, 2020.

\bibitem{koloskova2019decentralized}
A.~Koloskova, S.~Stich, and M.~Jaggi, ``Decentralized stochastic optimization
  and gossip algorithms with compressed communication,'' in \emph{International
  Conference on Machine Learning}.\hskip 1em plus 0.5em minus 0.4em\relax PMLR,
  2019, pp. 3478--3487.

\bibitem{nedic2008distributed}
A.~Nedi{\'c}, A.~Olshevsky, A.~Ozdaglar, and J.~N. Tsitsiklis, ``Distributed
  subgradient methods and quantization effects,'' in \emph{2008 47th IEEE
  Conference on Decision and Control}, 2008, pp. 4177--4184.

\bibitem{besbes2015non}
O.~Besbes, Y.~Gur, and A.~Zeevi, ``Non-stationary stochastic optimization,''
  \emph{Operations Research}, vol.~63, no.~5, pp. 1227--1244, 2015.

\bibitem{shalev2011online}
S.~Shalev-Shwartz \emph{et~al.}, ``Online learning and online convex
  optimization,'' \emph{Foundations and Trends in Machine Learning}, vol.~4,
  no.~2, pp. 107--194, 2011.

\bibitem{9184135}
X.~Li, X.~Yi, and L.~Xie, ``Distributed online optimization for multi-agent
  networks with coupled inequality constraints,'' \emph{IEEE Transactions on
  Automatic Control}, vol.~66, no.~8, pp. 3575--3591, 2021.

\bibitem{kalhan2021dynamic}
D.~S. Kalhan, A.~S. Bedi, A.~Koppel, K.~Rajawat, H.~Hassani, A.~K. Gupta, and
  A.~Banerjee, ``Dynamic online learning via {Frank-Wolfe} algorithm,''
  \emph{IEEE Transactions on Signal Processing}, vol.~69, pp. 932--947, 2021.

\end{thebibliography}

\end{document}